\documentclass[12pt]{amsart}
\theoremstyle{plain}
\newtheorem{theorem}{Theorem}[section]
\newtheorem{corollary}[theorem]{Corollary}
\newtheorem{example}[theorem]{Example}

\newtheorem{proposition}[theorem]{Proposition}
\newtheorem{lemma}[theorem]{Lemma}

\theoremstyle{definition}

\newtheorem{remark}[theorem]{Remark}
\newtheorem{problem}[theorem]{Problem}

\newcommand{\bb}[1]{\mathbb{#1}}
\newcommand{\cl}[1]{\mathcal{#1}}

\begin{document}

\title[Injective and projective Hilbert C*-modules]{Injective and projective
   Hilbert C*-modules, and C*-algebras of compact operators}
\author[M.~Frank]{Michael Frank}
\address{HTWK Leipzig, Fachbereich IMN, PF 301166, D-04251 Leipzig, F.R. Germany}
\email{mfrank@imn.htwk-leipzig.de}
\thanks{The authors were supported in part by a NSF grant.}
\author[V.~I.~Paulsen]{Vern I.~Paulsen}
\address{Dept.~Mathematics, Univ.~of Houston, Houston, TX 77204, U.S.A.}
\email{vern@math.uh.edu}
\begin{abstract}
We consider projectivity and injectivity of Hilbert C*-modules in
the categories of Hilbert C*-(bi-)modules over a fixed C*-algebra
of coefficients (and another fixed C*-algebra re\-presented as
bounded module operators) and bounded (bi-)module mor\-phisms,
either necessarily adjointable or arbitrary ones. As a consequence of these
investigations, we obtain a set of equivalent conditions
characterizing C*-subalgebras of C*-algebras of compact operators
on Hilbert spaces in terms of general properties of Hilbert
C*-modules over them. Our results complement results recently
obtained by B.~Magajna, J.~Schweizer and M.~Kusuda. In particular,
all Hilbert C*-(bi-)\-modules over C*-algebras of compact
operators on Hilbert spaces are both injective and projective in
the categories we consider. For more general C*-algebras we obtain
classes of injective and projective Hilbert C*-(bi-)modules.
\end{abstract}
\maketitle

The goal of this paper is to determine the injective and projective
Hilbert C*-modules over a fixed C*-algebra, when
one allows the maps between C*-modules to be bounded. Most prior work
on injectivity has focused on the case of contractive maps.  To better
understand our motivations and the distinction between this work and
the work of others, we first review the concept of {\em injectivity} 
and some history of the subject.

To give a definition of the term, {\em injective,} that is useful
for our purposes, we need a category, consisting of objects that
are sets and morphisms between them that are functions, and for
each object, $\cl N$, certain subsets, $\cl M \subseteq \cl N$
that are also objects in the category, which we call the {\em
subobjects} of $\cl N$. Then an object $\cl I$ in this category is
called {\em injective}, provided that for every object $\cl N$,
every subobject, $\cl M \subseteq \cl N$ and every morphism,
$\phi: \cl M \to \cl I$, there is a morphism $\psi: \cl N \to \cl
I$, that extends $\phi$. Note that if we keep the objects and
morphisms the same, but change the subobjects, then it is possible
that the injectives might change. Among the differences between
this definition of injectivity and the categorical definition that
is given in say, \cite{MacL}, is that our definition allows for
the possibility that the inclusion map of the subobject into the
object is not a morphism in the category. Many definitions of
injectivity, implicitly only consider subobjects with the property
that the inclusion map is a morphism.

If one fixes a ring $\cl R$, and considers the category whose
objects are left $\cl R$-modules, subobjects are left $\cl
R$-submodules and morphisms are left $\cl R$-module maps, then the
above concept of injective reduces to the classical definition of
an {\em injective left $\cl R$-module.} In this case the
inclusions of subobjects into objects are always morphisms.

Consider the category whose objects are Banach spaces, subobjects are
closed subspaces and whose morphisms are the contractive, linear
maps. Then it is easy to see, by a simple scaling, that a Banach
space, $\cl I$ is injective in this setting if and only if for
every Banach space, $\cl N$, every closed subspace, $\cl M
\subseteq \cl N$ and every bounded linear map, $\phi: \cl M \to
\cl I$, there is an extension, $\psi: \cl N \to \cl I$, of $\phi$
with $\|\phi\|= \|\psi\|$. A classic result, often called the
Nachbin-Goodner-Kelley theorem, is that a Banach space is
injective in this category if and only if it is isometrically
isomorphic to the space of continuous functions on an extremally
disconnected, compact Hausdorff space \cite{Ke}.

However, if one changes the category slightly, keeping the objects
to be Banach spaces and subobjects to be closed subspaces, but
allowing the morphisms to be all bounded, linear maps, then a
Banach space is injective in this category if and only if every
bounded, linear map on a subspace has a bounded, linear extension
to the whole space, but not necessarily of the same norm. A
complete understanding of the injective Banach spaces in this
setting is still unknown \cite{Ro}, but it is fairly easy to see
that Hilbert spaces are not injective.

Thus, generally, when one allows bounded maps to be the morphisms
instead of contractive maps, then the problems become more difficult.

On the other hand, if we now keep our morphisms to be bounded,
linear maps, but restrict our objects, by considering only Hilbert
spaces, then it is fairly easy to show that every Hilbert space
$\mathcal H$ is injective. This follows, since the extension can
be achieved by projecting onto the subspace. Thus, in this
restricted category, every object is now injective.

One way to generalize Hilbert spaces, is to consider the category
whose objects are Hilbert C*-modules over a fixed C*-algebra, $A,$
subobjects are Hilbert C*-submodules and morphisms are the
bounded, $A$-module maps. When $A= \bb C,$ then this category
reduces back down to the category of Hilbert spaces and bounded
linear maps.

Thus, in parallel with the Banach space case, we wish to determine
the injective objects in this setting. The first question that we
shall address is characterizing the C*-algebras $A$, such that,
like $\bb C$, every Hilbert $A$-module is injective.

As is the case with Banach spaces, if one restricts the morphisms
to be the contractive, module maps, then the theory of injective
Hilbert $A$-modules is somewhat simpler and is largely worked out
in the work of Huaxin Lin \cite{Lin1, Lin2} and Zhou Tian Xu
\cite{Xu}. In some of Lin's work, he studies injectivity, where in
our language, the objects are Hilbert C*-modules, subobjects are
Hilbert C*-submodules, and the morphisms are adjointable,
contractive module maps. In this case the inclusion map of a
submodule into the larger module is, generally, not a morphism,
since it need not be adjointable. Thus, to encompass the type of
``injectivity'' studied by Lin, one needs the more general
definition of injectivity given above and some care must be taken
when citing general facts about injectives from category theory in
his context. As we will show later, if the inclusion map is
required to be adjointable, i.e., if one restricts the subobjects,
then the submodules are necessarily orthogonally complemented and
every object is injective (Theorem \ref{basic-inj}).
Thus, the differences between our results and
those due to H.~Lin \cite{Lin1,Lin2} are caused by differences in
the categories that we consider.

So far we have only discussed injectivity, but similar comments
apply to {\em projectivity,} which in many ways is a dual theory
to injectivity. for the concept of projectivity, in addition to 
specifying the morphisms, one needs to specify the
{\em quotients.} We will make precise definitions of projectivity in
Section 4.

We shall also answer many parallel questions about characterizing
projective modules.

In the settings that we shall consider, the set of objects of all
cate\-gories under consideration consists of Hilbert C*-modules $\{
{\mathcal M}, \langle .,. \rangle \}$ over some fixed C*-algebra
$A$, i.e.~(left) $A$-modules $\mathcal M$ equipped with an
$A$-valued inner product $\langle .,. \rangle: {\mathcal M} \times
{\mathcal M} \to A$, cf.~\cite{Lance}. We specify a second
C*-algebra $B$ that is supposed to act on $\mathcal M$ as a set of
bounded adjointable operators via module-specific
$*$-representations. Thus, $\cl M$ is an $A$-$B$-bimodule, with
the right action of $B$ given by bounded adjointable maps on $\cl
M,$ so that, in particular, $\langle am_1b, m_2 \rangle= a \langle
m_1, m_2b^* \rangle,$ for every, $a \in A, b \in B,$ and $ m_1,m_2
\in \cl M.$ We call $\cl M$ a Hilbert $A$-$B$-bimodule.  The
requirement of the existence of a second action by $B$ changes the
unitary equivalence classes of Hilbert $A$-$B$-bimodules, i.e.,
the notion of equivalence in the categories under consideration.
Note that every Hilbert $A$-module is automatically a Hilbert
$A$-$\bb C$-module, where $\bb C$ denotes the complex numbers.

The sets of morphisms that we study will consist of either all
bounded bimodule morphisms between the objects, or all adjointable
bounded bimodule morphisms between them. We shall denote these two
categories by $\cl B(A,B)$ and $\cl B^*(A,B),$ respectively.

The subobjects that we will consider will be, generally, all
Hilbert $A$-$B$-submodules and occassionally the orthogonally
complemented $A$-$B$-submodules.

The primary goal of the present paper is the investigation of two
problems: $\,$ (i) characterize the C*-algebras $A$ and $B$ for
which any Hilbert $A$-$B$ bimodule is injective or projective for
one of the sets of morphisms under consideration and one of the
two concepts of subobjects; $\,$ (ii) find suitable sets of
projective, or injective, Hilbert $A$-$B$ bimodules for given
C*-algebras $A$ and $B$ and fixed morphism sets.

In most cases the action of the C*-algebra $B$ of bounded
$A$-linear operators on the Hilbert $A$-modules $\mathcal M$ turns
out to play a minor role. So we can concentrate on the C*-algebra
$A$ of module coefficients, on the Hilbert $A$-$B$ modules.
We obtain a full
characterization of the C*-algebras $A$ for which any Hilbert
$A$-$B$ bimodule is injective for both the principal
categories. For the category with only bounded adjointable $A$-$B$
bimodule maps as morphisms any C*-algebra of coefficients $A$ (and
any C*-algebra of bounded adjointable operators $B$) will suffice,
whereas for the category with all bounded $A$-$B$ bimodule maps as
morphisms only C*-algebras $A$ of compact operators (and arbitrary
C*-algebras $B$) have this property. If the C*-algebra of
coefficients $A$ is monotone complete then a Hilbert $A$-$B$
bimodule is injective in the category with the set of bounded
$A$-$B$ module maps if and only if it is self-dual.

In the case of projectivity, we show that every Hilbert $A$-$B$
bimodule is projective in the category $\cl B^*(A,B)$, for every
C*-algebra $A.$ We prove that when $A$ is a C*-algebra of compact
operators, then every Hilbert $A$-$B$ bimodule is projective in
$\cl B(A,B),$ but we are unable to resolve if these are the only
C*-algebras with this property. A characterization of such
algebras is not available at present. Even more, the question
whether all Hilbert C*-modules are projective in the larger
category, or not, remains open.

However, we do prove that all Hilbert $A$-$B$ bimodules over a
certain C*-algebra are projective in the categories investigated
 if and only if the kernel of every surjective bounded
module map between Hilbert $A$-modules is a topological direct
summand of the domain. Moreover, we identify a family of
projective C*-modules over unital C*-algebras. We show that the
finitely generated Hilbert C*-modules over unital C*-algebras
are projective objects of the categories under consideration.


There are some parallels between our results on projectivity and
research in progress on extensions of Hilbert C*-modules and on
projectivity of Hilbert C*-modules in this different category by
Damir Baki\'c and Boris Gulja{\v{s}}, \cite{BaGu1,BaGu2}.

Another way to modify the categories under consideration would be
to restrict the set of objects to self-dual (or orthogonally
comparable) Hilbert C*-modules. Recall that a Hilbert C*-module
$\cl M$ is {\em orthogonally comparable} provided that any time
$\phi: \cl M \to \cl N$ is an isometric module map, then $\phi(\cl
M)$ is orthogonally complemented in $\cl N$. However, this choice
implies the adjointability of all bounded module morphisms between
them and, consequently, that any Hilbert C*-submodule is an
orthogonal summand, cf.~\cite{Fr98}. So our questions would have
an immediate answer: in these latter categories all objects are
projective and injective for arbitrary C*-algebras $A$ and $B$.

Because of the close relation of the Magajna-Schweizer theorem
(\cite{Mag,Schweiz}) to
the circle of questions studied in the present paper, C*-algebras
$A$ of the form $A= c_0$-$\sum_\alpha \oplus {\rm K}(H_\alpha)$
are of special interest. Here the symbol ${\rm K}(H_\alpha)$
denotes the C*-algebra of all compact operators on some Hilbert
space $H_\alpha$, and the $c_0$-sum is either a finite block-diagonal
sum or a block-diagonal sum with a $c_0$-convergence condition on
the C*-algebra components ${\rm K}(H_\alpha)$. The $c_0$-sum may
possess arbitrary cardinality. These C*-algebras have been
precisely characterized by W.~Arveson \cite[\S I.4,
Th.~I.4.5]{Arveson} as C*-subalgebras of (full) C*-algebras of
compact operators on Hilbert spaces. We give a number of further
equivalent characterizations of this class of C*-algebras in terms
of general properties of Hilbert C*-modules over them which are of
separate interest. Throughout the present paper we refer to these
C*-algebras as {\it C*-algebras of compact operators on certain
Hilbert spaces}.

\section{Preliminaries}

In this section we give some definitions and basic facts from Hilbert
C*-module theory needed for our investigations. The papers
\cite{Pa1,Kas,Fr90,Lin:90/2,Lin1,Fr98}, some chapters in
\cite{JT,NEWO}, and the books by E.~C.~Lance \cite{Lance} and by
I.~Raeburn and D.~P.~Williams \cite{RaeWil} are used as standard
reference sources. We make the convention that all C*-modules
of the present paper are left modules by definition. A {\it
pre-Hilbert $A$-module over a C*-algebra} $A$ is an $A$-module
$\mathcal M$ equipped with an $A$-valued mapping $\langle .,.
\rangle : {\mathcal M} \times {\mathcal M} \rightarrow A$ which is
$A$-linear in the first argument and has the properties:
\[
\langle x,y \rangle = \langle y,x \rangle^* \; \, , \: \; \langle
x,x \rangle \geq 0 \quad {\rm with} \: {\rm equality} \: {\rm iff}
\quad x=0 \, .
\]
The mapping $\langle .,. \rangle$ is said to be {\it the
$A$-valued inner product on} $\mathcal M$. A pre-Hilbert
$A$-module $\{ \mathcal M, \langle .,. \rangle \}$ is {\it
Hilbert} if and only if it is complete with respect to the norm
$\| . \| = \| \langle .,. \rangle \|^{1/2}_A$. We always assume
that the linear structures of $A$ and $\mathcal M$ are compatible.
Two Hilbert $A$-modules are {\it isomorphic} if they are
isometrically isomorphic as Banach $A$-modules, if and only if
they are unitarily isomorphic, \cite{Lance}. We would like to
point out that an $A$-module can carry unitarily non-isomorphic
$A$-valued inner products which induce equivalent complete
norms, \cite{Fr98}. Two Hilbert $A$-$B$-modules are {\it isomorphic}
if and only if they are unitarily isomorphic as Hilbert
$A$-modules in such a way that the isomorphism intertwines the
$*$-representations of $B$ on them.

Hilbert C*-submodules of Hilbert C*-modules might not be direct
summands, and if they are direct summands then they might be
merely topological, but not orthogonal summands. We say that a
Hilbert C*-module $\mathcal N$ is a topological summand of a
Hilbert C*-module $\mathcal M$ which contains $\mathcal N$ as a
Banach C*-submodule in case $\mathcal M$ can be decomposed into
the direct sum of the Banach C*-submodule $\mathcal N$ and of
another Banach C*-submodule $\mathcal K$. The denotation is
${\mathcal M} = {\mathcal N} \stackrel{.}{+} {\mathcal K}$. If,
moreover, the decomposition can be arranged as an orthogonal one
(i.e.~$\mathcal N \bot \mathcal K$) then the Hilbert C*-submodule
$\mathcal N \subseteq \mathcal M$ is an orthogonal summand of the
Hilbert C*-module ${\mathcal M}$ i.e. ${\mathcal M} = {\mathcal K}
\oplus {\mathcal L}$.  Examples where these situations appear can
be found e.g.~in \cite{Fr98}.

A Hilbert $A$-module $\{ \mathcal M , \langle .,. \rangle \}$ over
a C*-algebra $A$ is said to be {\it self-dual} if and only if
every bounded module map $ r: \mathcal M \to A$ is of the form
$\langle . , x_r \rangle$ for some element $x_r \in \mathcal M$.
The set of all bounded module maps $r: \mathcal M \to A$ forms a
Banach $A$-module ${\mathcal M}'$. The module action of $A$ on
${\mathcal M}'$ is defined by the formula $(a \cdot r)(x) = r(x)a^*$
for any $x \in \mathcal M$, each $a \in A$ and $r \in {\mathcal M}'$.
A Hilbert $A$-module is called {\it C*-reflexive} (or more precisely,
{\it $A$-reflexive}) if and only if the map $\Omega$ defined by the
formula $\Omega (x)[r]=r(x)$ for each $x \in {\mathcal M}$, every
$r \in {\mathcal M}'$, is a surjective module map of $\mathcal M$
onto the Banach $A$-module ${\mathcal M}''$, where ${\mathcal M}''$
consists of all bounded module maps from ${\mathcal M}'$ to $A$.
Note that the property of being self-dual does not depend on the choice
of the C*-algebra of coefficients $A$ within $\langle {\mathcal M},
{\mathcal M} \rangle \subseteq A \subseteq M(\langle {\mathcal
M},{\mathcal M} \rangle)$, whereas the property of being
$A$-reflexive sometimes does.

As an example consider the
C*-algebra $A=c_0$ of all sequences converging to zero and set
$\mathcal M =c_0$ with the standard $A$-valued inner product.
Consider $\mathcal M$ both as a Hilbert $A$-module and as a
Hilbert $M(A)$-module. The multiplier C*-algebra of $A=c_0$ is
$M(A)=l_\infty$, the set of all bounded sequences. Then ${\mathcal M}'$
equals $l_\infty$ as a one-sided $A$-module, independently
of the choice of sets of coefficients. In contrast, the set of all
bounded $A$-linear maps of ${\mathcal M}'$ to $A$ can be identified
with $A=c_0$, whereas the set of all bounded $M(A)$-linear maps of
${\mathcal M}'$ to $M(A)$ can be identified with $l_\infty$. Generally
speaking, the $A$-dual Banach $A$-module ${\mathcal M}'$ of a
Hilbert $A$-module $\mathcal M$ can be described as the completion
of the linear hull of the unit ball of $\mathcal M$ with respect to the
topology induced by the seminorms, $\{ \|\langle x,. \rangle \|_A
: x \in {\mathcal M}, \, \|x\| \leq 1 \}$ \cite[Th.~6.4]{Fr98}.
The process of forming higher order C*-dual Banach C*-modules of a
given Hilbert C*-module $\mathcal M$ stabilizes after the second
step since ${\mathcal M}' \equiv {\mathcal M}'''$. We have the
standard chain of isometric Banach C*-module embeddings ${\mathcal
M} \subseteq {\mathcal M}'' \subseteq {\mathcal M}'$ by
\cite{Pa1,Pa2}.

Furthermore, we are going to consider various bounded C*-linear
operators $T$ between Hilbert C*-modules $\mathcal M$, $\mathcal
N$ with one and the same C*-algebra of coefficients. Quite
regularly those operators $T$ may not admit an adjoint bounded
C*-linear operator $T^*: \mathcal N \to \mathcal M$ fulfilling the
equality $\langle T(x),y \rangle_{\mathcal N} = \langle x, T^*(y)
\rangle_{\mathcal M}$ for any $x \in \mathcal M$, any $y \in
\mathcal N$. We denote the C*-algebra of all bounded C*-linear
adjointable operators on a given Hilbert $A$-module $\mathcal M$
by ${\rm End}_A^*(\mathcal M)$. The Banach algebra of all bounded
$A$-linear operators on $\mathcal M$ is denoted by ${\rm
End}_A(\mathcal M)$. For more detailed information on these
situations we refer to \cite{Fr98}.

A result that we shall use often is a bounded closed graph theorem
for Hilbert C*-modules that is a variant of N.~E.~Wegge-Olsen's
result. We show how the bounded closed graph theorem can be
derived from his result. In contrast, an example of E.C. Lance
shows that there is no analogue of the unbounded closed graph
theorem for general Hilbert C*-modules.

\begin{proposition} \label{closed-range} {\rm (N.~E.~Wegge-Olsen
                          \cite[Th.~15.3.8]{NEWO})} \newline
   Let $A$ be a C*-algebra, $\{ {\mathcal M}, \langle .,. \rangle \}$
   be a Hilbert $A$-module and $T$ be an adjointable bounded module
   operator on $\mathcal M$.
   If $T$ has closed range then $T^*$, $(T^*T)^{1/2}$ and $(TT^*)^{1/2}$
   have also closed ranges and
   \begin{eqnarray*}
   {\mathcal M} & = & {\rm Ker}(T) \oplus T^*({\mathcal M})
              =   {\rm Ker}(T^*) \oplus T({\mathcal M}) \\
            & = & {\rm Ker}(|T|) \oplus |T|({\mathcal M})
              =   {\rm Ker}(|T^*|) \oplus |T^*|({\mathcal M}) \, .
   \end{eqnarray*}
   In particular, each orthogonal summand appearing on the right is
   automatically norm-closed and coincides with its bi-orthogonal
   complement inside $\mathcal M$. Moreover, $T$ and $T^*$ have polar
   decomposition.
\end{proposition}

\begin{corollary}  {\rm (bounded closed graph theorem)} \label{adjointability}
    \newline
   Let $A$ be a C*-algebra and $\{ \mathcal M, \langle .,. \rangle \}$,
   $\{ \mathcal N, \langle .,. \rangle \}$ be two Hilbert $A$-modules.
   The graph of every bounded
   $A$-linear operator $T$ coincides with its bi-orthogonal complement in
   $\mathcal M \oplus \mathcal N$, and it is always a topological summand
   with topological complement $\{ (0,z) : z \in \mathcal N \}$.
   A
   bounded $A$-linear operator $T: \mathcal M \to \mathcal N$ possesses an
   adjoint operator $T^*: \mathcal N \to \mathcal M$ if and only if the
   graph of $T$ is an orthogonal summand of the Hilbert $A$-module
   $\mathcal M \oplus {\mathcal N}$.
\end{corollary}

\begin{remark}
   By a counterexample due to E.~C.~Lance (\cite[pp.~102-104]{Lance}) the
   graph of a closed, self-adjoint, densely defined, unbounded module
   operator need not coincide with its bi-orthogonal complement.
\end{remark}

\begin{proof}
Since the inequality $\|T(x)\| \leq \|T\| \|x\|$ is valid for every $x \in
\mathcal M$ the graph of $T$ is a norm-closed Hilbert $A$-submodule of the
Hilbert $A$-module $\mathcal M \oplus \mathcal N$. Moreover, since the graph
of $T$ is the kernel of the bounded module operator $S: (x,y) \to (0,T(x)-y)$
on $\mathcal M \oplus \mathcal N$ it coincides with its bi-orthogonal
complement there, \cite[Cor.~2.7.2]{Fr97}. If $T$ has an adjoint then the
operator $T':(x,y) \to (x,T(x))$ is adjointable on $\mathcal M \oplus
\mathcal N$. By Proposition \ref{closed-range} the graph of $T$ is an
orthogonal summand.

Conversely, if the graph of $T$ is an orthogonal summand of
${\mathcal M} \oplus {\mathcal N}$ then its orthogonal complement
consists precisely of the pairs of elements $\{ (x,y) : x =
-T^*(y), y \in {\mathcal N} \}$. To see this consider the equality
$\langle z,x \rangle_{\mathcal M} + \langle T(z),y
\rangle_{\mathcal N}=0$ which has to be valid for any $z \in
\mathcal M$ and any pair $(x,y)$ of elements of the orthogonal
complement of the graph of $T$. The assumption of the existence of
two pairs $(x_1,y)$ and $(x_2,y)$ in this complement forces
$\langle z, x_1 \rangle_{\mathcal M} = \langle z, x_2
\rangle_{\mathcal M}$ for any $z \in \mathcal M$, and therefore,
$x_1=x_2$. Hence, $T^*(y)(z) = \langle z,(-x) \rangle_{\mathcal
M}$ for any $z \in \mathcal M$ and for $T^*:{\mathcal N} \to
{\mathcal M}'$. So $T^*$ is everywhere defined on $\mathcal N$
taking values exclusively in ${\mathcal M} \subseteq {\mathcal
M}'$. This shows the existence of the adjoint operator $T^*$ of
$T$ in the sense of its definition.

The property of the graph of a bounded module operator to be a topological
summand with topological complement $\{ (0,z) : z \in \mathcal N \}$ follows
from the decomposition $(x,y) = (x,T(x)) + (0,y-T(x))$ for every $x \in
\mathcal M$, $y \in \mathcal N$. Since $T(0)=0$ for any linear operator $T$
the intersection of the graph with the $A$-$B$-submodule $\{ (0,z) : z \in
\mathcal N \}$ is always trivial.
\end{proof}

There is still one open problem about complements whose solution for (at
least, surjective) bounded module mappings would give us insight
into the solution of the main question of the
fourth section concerning projective Hilbert C*-modules.

\begin{problem}
 Suppose, a bounded module operator between Hilbert C*-modules
 has a norm-closed image which is either a topological summand
 or merely coincides with its biorthogonal complement with respect
 to the image Hilbert C*-module.
 Is the kernel of such an operator always a topological summand,
 or are there counterexamples?
\end{problem}

The difficulties surrounding this problem are illuminated by an example
constructed by V.~M.~Manuilov in \cite{Manuilov:00}.

\section{C*-algebras of compact operators and the Magajna-Schweizer theorem}

In this section, we prove that
the class of C*-algebras of compact operators on certain Hilbert spaces and
their C*-subalgebras can be characterized by the appearance of certain properties
common to all Hilbert C*-modules over them. The different aspects
shown below enable us to establish classes of injective
and projective Hilbert $A$-modules for these and other C*-algebras of
coefficients $A$ in forthcoming sections. Our starting point is the following
result by Bojan Magajna and J\"urgen Schweizer:

\begin{theorem} \label{Mag-Schw} {\rm (B.~Magajna, J.~Schweizer
                     \cite{Mag,Schweiz})} \newline
   Let $A$ be a C*-algebra. The following three conditions are equivalent:

   \newcounter{cou001}
   \begin{list}{(\roman{cou001})}{\usecounter{cou001}}
   \item  $A$ is of $c_0$-$\sum_i \oplus {\rm K}(H_i)$-type, i.e.~it has a
          faithful $*$-representation as a C*-algebra of compact operators
          on some Hilbert space.
   \item  For every Hilbert $A$-module $\mathcal M$ every Hilbert $A$-submodule
          \linebreak[4]
          $\mathcal N \subseteq \mathcal M$ is automatically orthogonally
          complemented in $\mathcal M$, i.e.~$\mathcal N$ is an orthogonal
          summand of $\mathcal M$.
   \item  For every Hilbert $A$-module $\mathcal M$ every Hilbert $A$-submodule
          \linebreak[4]
          $\mathcal N \subseteq \mathcal M$ that coincides with its
          bi-orthogonal complement ${\mathcal N}^{\bot\bot} \subseteq \mathcal
          M$ is automatically orthogonally complemented in $\mathcal M$.
   \end{list}
\end{theorem}

Based on the Magajna-Schweizer theorem further investigations were made for
the identification of generic general properties of Hilbert C*-modules which
characterize entire classes of C*-algebras of coefficients, cf.~\cite{Fr05}.
Many of these generic properties turned out to characterize C*-algebras of
compact operators in case they are common for all Hilbert C*-modules over a
certain C*-algebra of coefficients. We present these properties here as a
list of equivalent conditions that extend the conditions of Magajna-Schweizer.

\begin{proposition} \label{Fr-neu} \cite{Fr05}
   Let $A$ be a C*-algebra. The following seven conditions are equivalent:

   \begin{list}{(\roman{cou001})}{\usecounter{cou001}}
   \item  $A$ is of $c_0$-$\sum_i \oplus {\rm K}(H_i)$-type, i.e., it has a
          faithful $*$-representation as a  C*-algebra of compact operators
          on some Hilbert space.
   \setcounter{cou001}{3}
   \item For every Hilbert $A$-module $\mathcal M$ and every bounded $A$-linear
         map $T: {\mathcal M} \to {\mathcal M}$ there exists an adjoint bounded
         $A$-linear map $T^*:{\mathcal M} \to {\mathcal M}$.
   \item For every pair of Hilbert $A$-modules $\mathcal M$, $\mathcal N$ and
         every bounded $A$-linear map $T: {\mathcal M} \to {\mathcal N}$ there
         exists an adjoint bounded $A$-linear map $T^*:{\mathcal N} \to
         {\mathcal M}$.
   \item The kernels of all bounded $A$-linear operators between arbitrary
         Hilbert $A$-mo\-du\-les are orthogonal summands.
   \item The images of all bounded $A$-linear operators with norm-closed
         range between arbitrary Hilbert $A$-modules are orthogonal summands.
   \item For every Hilbert $A$-module every Hilbert $A$-submodule
         is automatically topologically complemented, i.e., it is a
         topological summand.
   \item For every (maximal) norm-closed left ideal $I$ of $A$ the
         corresponding open projection $p \in A^{**}$ is an element of the
         multiplier C*-algebra ${\rm M}(A)$ of $A$.
   \end{list}
\end{proposition}

We will see in the following sections that some of these
equivalent conditions force Hilbert C*-modules over C*-algebras of
compact operators to be projective or injective. The investigation
of these generic categorical properties of Hilbert C*-modules
revealed, however, a problem that is still unsolved. It is related
to the identification of (non-) injective and (non-)projective
Hilbert C*-modules, and so we list it here:

\begin{problem} \label{problemYY}
 Characterize those C*-algebras $A$ for which the following condition holds:
 For every Hilbert $A$-module, every Hilbert $A$-submodule that
 coincides with its bi-ortho\-gonal complement is automatically
 topologically complemented there.
\end{problem}

Problem \ref{problemYY} revisits the difference between
B.~Magajna's theorem and J.~Schwei\-zer's theorem on the level of
topological summands. The results by M. Kusuda \cite{Kus2}
indicate that the solutions of these problems has to be similar to
his results on orthogonal summands. M.~Kusuda considered the
problem in \cite{Kus2} in 2005 and has got a number of results
towards a solution in the spirit of the Magajna-Schweizer theorem
using spectral methods for C*-algebras. However, his results
indicate that the final solution of Problem \ref{problemYY} might
not have a simple formulation but might consist of a rather
extended list of cases to be distinguished.

\section{Injectivity}

\noindent
Let $A$ and $B$ be two fixed C*-algebras. We consider
two categories. In both categories the {\em objects} will be the
Hilbert $A$-$B$-bimodules. The sets of {\em morphisms} that we
study will consist of either all bounded bimodule morphisms
between the objects, or all adjointable, bounded bimodule
morphisms between them. In both cases, the {\em subobjects}  will
be the set of all Hilbert $A$-$B$-submodules, that is, norm closed
subspaces which are invariant under both the module actions.

We shall denote these two categories, together with the specified
sets of subobjects, by $\mathcal B(A,B)$ and $\mathcal B^*(A,B)$,
respectively. Note that every left $A$-module is always equipped
with a (right) action by the complex numbers, $\mathbb C$. Thus, $\mathcal
B(A, \mathbb C)$ (respectively, $\mathcal B^*(A, \mathbb C)$) is just the
category of left Hilbert $A$-modules and bounded (respectively,
bounded, adjointable) maps.

So, in summary, a Hilbert $A$-$B$-bimodule $\mathcal E$ is {\it
injective in} $\mathcal B(A,B)$, (respectively, $\mathcal B^*(A,B)$) if and
only if for every Hilbert $A$-$B$-bimodule, $\mathcal N$, and every
Hilbert $A$-$B$-subbimodule, $\mathcal M$ of $\mathcal N$, and every
bounded, (respectively, bounded, adjointable) bimodule map, $\phi:
\mathcal M \to \mathcal E$, there is a bounded (respectively, bounded,
adjointable), bimodule map $\psi: \mathcal N \to \mathcal E$
that extends $\phi$. In other words, a Hilbert $A$-$B$ bimodule
$\mathcal E$ is injective if and only if the diagram
\begin{equation}  \label{diag-inj}
   \begin{array}{lcc}
       {\mathcal N} & & \\
       \uparrow\vcenter{\rlap{$T$}} & & \\
       {\mathcal M} & \stackrel{\phi}{\rightarrow} & {\mathcal E}
   \end{array}
\end{equation}
can be completed to a commutative one by an $A$-$B$ bimodule morphism $\psi:
\mathcal N \to \mathcal E$ of the selected category.

Before beginning our study of injectivity, we first point out what
happens when the set of subobjects is required to be smaller.

The following theorem should be contrasted with H.~Lin's result
\cite[Th.~2.14]{Lin1}, which applies to the category of left
Hilbert $A$-modules with morphisms the contractive adjointable
maps, but a larger family of subobjects was allowed, namely, all
the submodules were considered subobjects. Thus, the inclusion
maps were not in general morphisms. H.~Lin obtained that, in this
setting, a Hilbert $A$-module is injective if and only if it is
orthogonally comparable as a Hilbert $A$-module. We can
demonstrate that even expanding the morphisms to the bounded
adjointable $A$-$B$-bimodule maps, but requiring the inclusion
maps to be morphisms, changes the picture rather significantly.

\begin{theorem} \label{basic-inj}
   Let $A$ be an arbitrary C*-algebra and $\{ {\mathcal E},
   \langle .,. \rangle \}$ be a Hilbert $A$-module.  Let $B$
   be another C*-algebra admitting a $*$-repre\-sen\-tation in
   ${\rm End}_A^*({\mathcal E})$. Then $\mathcal E$ is an
   injective object in the category whose objects are the Hilbert
   $A$-$B$-bimodules, whose morphisms are either the (adjointable)
   contractive or (adjointable) bounded bimodule maps and whose
   subobjects are the $A$-$B$-submodules whose inclusion maps are
   adjointable. Consequently, every element of those categories
   is injective.
\end{theorem}

\begin{proof}
Since by assumption the inclusion $T: {\mathcal M} \hookrightarrow
\mathcal N$ is an adjointable bounded $A$-$B$-bimodule map the map
$T^*$ is a surjective bounded $A$-$B$-bimodule map and Proposition
\ref{closed-range} applies: the image set $T({\mathcal M})
\subseteq \mathcal N$ is an orthogonal summand. Moreover, the map
$T^{-1}: T({\mathcal M}) \to \mathcal M$ defined as
$T^{-1}(T(x))=x$ for $x \in \mathcal M$ is everywhere defined on
$T({\mathcal M}) \subseteq \mathcal N$ and bijective, so it is
bounded and $A$-$B$-bilinear by definition. It can be extended to
a map defined on $\mathcal N$ simply setting it to be the zero map
on the orthogonal complement of $T( {\mathcal M})$ in $\mathcal
N$. Preserving the denotation $T^{-1}$ for this extension, setting
$\psi=\phi \circ T^{-1}$ yields the desired extension of $\phi$ to
$\mathcal N$. Consequently, the Hilbert $A$-$B$-bimodule $\mathcal E$
is automatically injective in the category under consideration.
\end{proof}

We now focus on the two categories that are our principal
interest.

To make further progress in identifying the injective objects of
the category ${\mathcal B}(A,B)$ we consider consequences of the
definition of injectivity.

\begin{lemma}   \label{lem-topsumm}
   Let $A$, $B$ be C*-algebras and $\{ {\mathcal E}, \langle .,.
   \rangle \}$ be an injective Hilbert $A$-$B$-bimodule in one of
   the two categories under consideration. If $\mathcal E \subseteq
   \mathcal N$ is an $A$-$B$-submodule, then the Hilbert $A$-$B$-bimodule
   $\mathcal E$ is a topological summand of the Hilbert $A$-$B$
   bimodule $\mathcal N$.
   \newline
   Moreover, $\mathcal E$ is $A$-reflexive as a Hilbert $A$-module,
   and whenever $\mathcal E$ is a Hilbert $A$-submodule of another
   Hilbert $A$-module $\mathcal M$ with ${\mathcal E}^\bot = \{ 0
   \}$ then $\mathcal E = {\mathcal E}^{\bot\bot}$ in $\mathcal M$.
\end{lemma}

\begin{proof} In the definition of injectivity, let $\mathcal M = \mathcal
E$, let $T$ denote the inclusion of $\mathcal E$ into $\mathcal N$ and let
$\phi={\rm id}_{\mathcal E}$. By supposition there exists an $A$-$B$
bimodule morphism $\psi: \mathcal N \to \mathcal E$ such that
$\psi \circ T = {\rm id}_{\mathcal E}$. By \cite[Lemma
3.1.8(2)]{Kasch} we have the set identities $\mathcal N =
\psi^{-1} (\mathcal E) = {\rm Im}(T) + {\rm Ker}(\psi)$ and $\{ 0
\} = T({\rm Ker} ({\rm id}_{\mathcal E})) = {\rm Im}(T) \cap {\rm
Ker}(\psi)$. Therefore, $\mathcal E$ has to be a topological
summand with topological complement ${\rm Ker}(\psi)$ there,
i.e.~$\mathcal N = T(\mathcal E) \stackrel{.}{+} {\rm Ker}(\psi)$.

To derive the $A$-reflexivity of injective Hilbert $A$-modules
consider the definition of injectivity with $\mathcal M = \mathcal
E$, $\mathcal N = {\mathcal E}''$ and $\phi = {\rm id}_{\mathcal
E}$. By \cite[Prop.~2.1]{Pa2} the $A$-valued inner product on
$\mathcal E$ extends to an $A$-valued inner product on its
$A$-bidual Banach $A$-module ${\mathcal E}''$. Moreover, the
$*$-representation of $B$ on $\mathcal E$ turns into a
$*$-representation of $B$ on ${\mathcal E}''$ via the canonical
isometric embedding $\mathcal E \subseteq {\mathcal E}''$, since
every bounded module operator on $\mathcal E$ extends to a bounded
module operator on ${\mathcal E}''$ in a unique way by \cite{Pa2}.
However, the embedded copy of $\mathcal E$ is a topological
summand of ${\mathcal E}''$ if and only if both they coincide.
Indeed, since we have the chain of isometric embeddings $\mathcal
E \subseteq {\mathcal E}'' \subseteq {\mathcal E}'$ by
\cite{Pa1,Pa2} the assumption of $\mathcal E$ being a non-trivial
topological summand of ${\mathcal E}''$ would lead to the
non-uniqueness of the representation of the zero map on $\mathcal
E$ in ${\mathcal E}'$, a contradiction to the definition of this
set. The last statement above is a consequence of the injectivity
and $A$-reflexivity of $\mathcal E$ and of \cite[Lemma
3.1]{Fr98/2}.
\end{proof}

In many categories the converse of Lemma \ref{lem-topsumm} holds
too, that is, if an object is complemented in every object that it
is a subobject of, then it is injective. This holds any time that
there are enough injective elements in the category that the
object can be embedded via a morphism as a subobject of an
injective object. Often these splitting properties serve as an
alternative means to define injectivity, \cite{Kasch}. However, in
Proposition \ref{prop-example} below we indicate that for a large
number of unital, monotone incomplete C*-algebras $A$ the
C*-algebra $A$ itself is not injective in the category of all
Hilbert $A$-modules and bounded module maps, despite the fact that
unital C*-algebras $A$ are always orthogonally comparable as
Hilbert $A$-modules. The same holds for certain non-unital
C*-algebras $A$ provided ${\rm M}(A) = {\rm LM}(A)$

\begin{proposition} \label{prop-example}
   Let $A$ be a C*-algebra and $A^N$ be the standard Hilbert
   $A$-module of all $N$-tuples of elements of $A$ for given
   positive integers $N$. The following are equivalent:
   \begin{list}{(\roman{cou001})}{\usecounter{cou001}}
   \item  $A^N$ is injective in $\mathcal B(A, \mathbb C)$ for one
          $N \in \mathbb N$,
   \item  $A^N$ is injective in $\mathcal B(A, \mathbb C)$ for every
          $N \in \mathbb N$,
   \item  $A$ is injective in $\mathcal B(A, \mathbb C)$,
   \item  ${\rm M}(A)$ is a monotone complete C*-algebra.
   \end{list}
\end{proposition}

\begin{proof}
Let $\mathcal M \subseteq \mathcal N$ be a subobject and let $\phi: \mathcal M
\to A^N$ be a bounded $A$-module map. Note that $\phi= (\phi_1,
\ldots, \phi_N)$ where $\phi_i: \mathcal M \to A$ are bounded
$A$-module maps. The map $\psi: \mathcal N \to A^N$ that extends $\phi$
exists if and only if there exist bounded $A$-module maps $\psi_i:
{\mathcal N} \to A_{(i)}$ coinciding with $\phi_i$ on $\mathcal
M$, where the index $(i)$ denotes $i$-th coordinate of $A^N$. This
shows the equivalence of (1), (2) and (3).

Also, we see that such an extension exists if and only if a
generalized Hahn-Banach type theorem is valid for arbitrary pairs
of Hilbert $A$-modules $\mathcal M \subseteq \mathcal N$ and
arbitrary bounded $A$-linear functionals $r: \mathcal M \to A$. By
\cite[Th.~2]{Fr98/2} this takes place if and only if ${\rm M}(A)$
is monotone complete.
\end{proof}

\begin{proposition} \label{prop-full}
  Let $A$ be a unital C*-algebra.
  If there exist any full Hilbert $A$-modules that are injective in
  $\mathcal B(A, \mathbb C)$, then $A$ is monotone complete. Hence,
  if $A$ is simple, unital and not monotone complete, then there are
  no non-zero injective Hilbert $A$-modules in $\mathcal B(A,
  \mathbb C)$.
\end{proposition}

\begin{proof}
Let $\mathcal E$ be a full injective Hilbert $A$-module. By
\cite[Lemma 2.4.3]{ManTroi} there exists a finite positive
integer $n$ and a subset of elements $\{ e_1, ... , e_n \}$
of $E$ such that $\sum_{i=1}^n \langle e_i,e_i \rangle = 1_A$
because the Hilbert $A$-module $E$ is full. Note, that $E^n$ is
injective whenever $E$ is injective and $n$ is a finite
positive integer. Therefore,
we have an isometric left $A$-module map, $\phi:A \to \mathcal E^n$
given by $\phi(a) = \sum_{i=1}^n ae_i$. Since, $A$ is orthogonally
comparable by \cite[Prop.~6.2, Th.~6.3]{Fr98}, there exists a bounded
$A$-module map, $\psi: \mathcal E^n \to \phi(A)$. From this it follows
easily that $\phi(A)$ is injective in $\mathcal B (A, \mathbb C)$.
Since $A$ and $\phi(A)$ are isomorphic, $A$ is injective in
$\mathcal B (A, \mathbb C)$. Hence, by Proposition \ref{prop-example},
${\rm M}(A)= A$ is monotone complete.

To see the final assertion, note that since $A$ is simple and unital,
every non-zero Hilbert $A$-module is full, since the range of its'
$A$-valued inner product is a norm-closed two-sided ideal in $A$.
\end{proof}

When $A$ is unital, not simple and not monotone complete, it is
possible to have injectives in $\mathcal B(A, \mathbb C)$, as the following
example shows. However, we will show below that when $A$ is unital
but not monotone complete, then there are not enough injectives, so
that every Hilbert $A$-module can be embedded in an injective.

\begin{example} {\rm Let $A= \mathbb C \oplus B$, where $B$ is a unital
  C*-algebra that is assumed to be not monotone complete. Thus, $A$ is unital
  and not monotone complete. Note, that every Hilbert space $\mathcal K$ is a
  (non-full) Hilbert $A$-module, with $(0 \oplus B) \mathcal K = 0$. We
  claim that $\mathcal K$ is an injective Hilbert $A$-module in $\mathcal B(A,
  \mathbb C)$.

  Indeed, in case $\mathcal E$ is a Hilbert $A$-module and $\mathcal H=
  (\mathbb C \oplus 0)\mathcal E$ and $\mathcal F= (0 \oplus B)\mathcal E$
  are its submodules, then $\mathcal E= \mathcal H \oplus \mathcal F$ is an
  orthogonal direct sum decomposition. Moreover, any $A$-module map from
  $\mathcal E$ into $\mathcal K$ is zero on $\mathcal F$, and it is a linear
  map on $\mathcal H$. The fact that $\mathcal K$ is injective in $\mathcal B
  (A, \mathbb C)$ now follows easily from that it is injective in the
  category of Hilbert spaces and bounded linear maps.}
\end{example}

By the above results, the category ${\mathcal B} (A,\mathbb C)$ does
not contain any non-zero injective object for setting $A$ to be one
of the following C*-algebras, among others (cf.~\cite{Davidson}):
\begin{itemize}
  \item{the reduced group C*-algebra $C_r^*(F_2)$ of the free group
        on two generators $F_2$,}
  \item{the irrational rotation algebras $A_\theta$, $\theta \in (0,1)$
        - irrational,}
  \item{the Cuntz algebras $O_n$, $n \in \mathbb N$ and the Cuntz-Krieger
        algebra $O_\infty$,}
  \item{the Bunce-Deddens algebras ${\mathcal B}(\{n_k\})$, $n,k \in
        \mathbb N$,}
  \item{Blackadar's projectionless unital simple C*-algebra.}
\end{itemize}

The following result, shows that even in cases when injectives do exist
for $A$ unital and not monotone complete, there cannot exist
``enough''. We say that a Hilbert $A$-module, $\mathcal E$ can be {\em
boundedly embedded} in a Hilbert $A$-module $\mathcal F$, provided that
there exists a module map, $T: \mathcal E \to \mathcal F$ that is bounded
above and below, i.e., there are constants, $0<C_1 \le C_2$ such
that $C_1 \|e\| \le \|T(e)\| \le C_2\|e\|$.

Often a category is said to have {\em enough injectives} if every
object can be embedded into an injective object. The following
result shows that if $A$ is unital, then the only time that $\mathcal
B(A, \mathbb C)$ can have enough injectives is when $A$ is monotone
complete.

\begin{proposition} \label{prop-enough}
  Let $A$ be a unital C*-algebra. If every Hilbert $A$-module can
  be boundedly embedded into a Hilbert $A$-module that is injective
  in $\mathcal B(A, \mathbb C)$, then $A$ is monotone complete.
\end{proposition}

\begin{proof}  By hypothesis, there exists an injective Hilbert
$A$-module $\mathcal F$ and a bounded embedding $T:A \to \mathcal F$
with constants $C_1, C_2$ as above. Since $T$ is an $A$-module map,
there exists an element $f \in \mathcal F$ such that $T(a) = af$.

Let $p= \langle f,f \rangle^{1/2}$. By the
above inequalities, $C_1^2 \|a\|^2 \le \|ap^2a^*\| \le C_2^2
\|a\|^2$ for every $a \in A$. Taking $a= g(p)$ for $g$ some
continuous function on the spectrum of $p$ and using the fact that
$p$ is a positive element of $A$, we get by the spectral mapping
theorem, that $C_1 \|g(t)\|_{\infty} \le \|tg(t)\|_{\infty} \le
C_2\|g(t)\|_{\infty}$, where $\|\cdot\|_{\infty}$ denotes the
supremum over the spectrum of $p$. These inequalities imply that
the spectrum of $p$ is contained in the interval, $[C_1, C_2]$,
and hence $p$ is invertible in $A$.

Therefore, $\mathcal F$ is an injective full Hilbert $A$-module
and by Proposition \ref{prop-full}, $A$ is monotone complete.
\end{proof}

\begin{problem} We do not have a complete set of analogous results
  for C*-bimodules. For example, if $A$ and $B$ are unital, simple
  C*-algebras, and neither one is monotone complete, can there exist
  any full Hilbert $A$-$B$-bimodules that are injective in $\mathcal B(A,B)$?
  Under what conditions on unital C*-algebras $A$ and $B$,
  can there exist enough injectives in $\mathcal B(A,B)$?
\end{problem}

In contrast, for non-unital C*-algebras there are often enough
injectives.

The following result characterizes the C*-algebras for which every
Hilbert $A$-module is injective in $\mathcal B(A, \mathbb C)$.

\begin{theorem} \label{inj_nonunit}
   Let $A$ be a C*-algebra of compact operators on some Hilbert
   space. Let $\{ {\mathcal E}, \langle .,. \rangle \}$ be a Hilbert
   $A$-module and $B$ be another C*-algebra admitting a
   $*$-representation on $\mathcal E$. Then $\mathcal E$ is an
   injective object in $\mathcal B(A, B)$.
   \newline
   Conversely, let $A$ be a C*-algebra. If every Hilbert $A$-module
   is injective in $\mathcal B(A,\mathbb C)$, then $A$ is
   $*$-isomorphic to a C*-algebra of compact operators on some
   Hilbert space.
\end{theorem}

\begin{proof}
Referring to Theorem \ref{Mag-Schw} and Proposition \ref{Fr-neu}
we see that every bounded $A$-linear map between Hilbert
$A$-modules over a C*-algebra $A$ of type $c_0$-$\sum_i \oplus
{\rm K}(H_i)$ possesses an adjoint. So every inclusion map is
adjointable and we are in the situation of Theorem \ref{basic-inj}.
This shows the first assertion.

To demonstrate the converse implication consider a maximal
left-sided ideal $I$ of the C*-algebra $A$. The $A$-valued inner
product on $I$ is that one inherited from $A$. Setting ${\mathcal E} =
{\mathcal M} = I$, ${\mathcal N} = A$, $\phi={\rm id}_I$ and
taking the standard $A$-linear embedding of $I$ into $A$ in
the definition of injectivity, we see that the existence of an
$A$-module map $\psi: A \to I$ extending $\phi$ is equivalent to
the existence of an orthogonal projection $p_I \in {\rm M}(A)$
such that $I = A p_I$. So by Proposition \ref{Fr-neu}, (ix) the
C*-algebra $A$ has to be of type $c_0$-$\sum_i \oplus {\rm K}(H_i)$
as to be shown.
\end{proof}

Let us remark that maximal left-sided ideals $I$ of C*-algebras $A$
may admit trivial sets of bounded module operators ${\rm End}_A(I)
= \mathbb C$. For example, consider the case of $A$ being the
matrix algebra of complex $2 \times 2$ matrices. Consequently, the
validity of the second assertion for the more general setting of
Hilbert $A$-$B$-bimodules heavily depends on the structure of the
C*-algebra $B$ and of the resulting collection of objects of
$\mathcal B (A, B)$.

When $A$ is a monotone complete
C*-algebra we can characterize the injectivity of Hilbert $A$-$B$
bimodules in terms of self-duality. This strengthens a result of
H.~Lin (\cite[Th.~2.2,Prop.~3.10]{Lin1}) that he obtained in the
contractive morphism situation, since every module that is injective
in the contractive morphism situation is automatically injective in
the setting of bounded morphisms. Also, our result complements a
result by D.~P.~Blecher and V.~I.~Paulsen on the injective
envelope of an operator bimodule stating that it has to be a
self-dual Hilbert C*-module over an injective (and hence, monotone
complete) C*-algebra, cf.~\cite{BP}.

\begin{theorem} \label{thm-criterion}
   Let $A$ be a monotone complete C*-algebra, $\{ \mathcal M,\langle
   .,. \rangle \}$ be a Hilbert $A$-module. Let $B$ be a C*-algebra
   admitting a $*$-representation in ${\rm End}_A^*({\mathcal M})$.
   Then $\mathcal M$ is injective in $\mathcal B(A,B)$ if and only
   if $\mathcal M$ is self-dual as a Hilbert $A$-module.
\end{theorem}

\begin{proof}
Suppose $\mathcal M$ is injective in $\mathcal B(A,B)$, and
consider the canonical isometric embedding of $\mathcal M$ into
its $A$-dual Banach $A$-module ${\mathcal M}'$. By
\cite[Th.~4.7]{Fr3} the $A$-valued inner product on $\mathcal M$
can be continued to an $A$-valued inner product on ${\mathcal M}'$
in a manner compatible with the canonical embedding ${\mathcal M}
\hookrightarrow {\mathcal M}'$. For the action on the right, the
$*$-representation of $B$ on $\mathcal M$ induces a
$*$-representation of $B$ on ${\mathcal M}'$ via the canonical
embedding, since every bounded module operator on $\mathcal M$
extends to a unique bounded module operator on ${\mathcal M}'$ by
\cite{Pa1}. Finally, the copy of $\mathcal M$ in ${\mathcal M}'$
is a topological summand there if and only if both the sets
coincide, since the zero functional on $\mathcal M$ would admit
several representations in ${\mathcal M}'$ otherwise. So $\mathcal
M$ has to be self-dual.

To establish the converse implication consider the diagram
\[
   \begin{array}{lcc}
       {\mathcal L} & & \\
       \uparrow\vcenter{\rlap{$T$}} & & \\
       {\mathcal K} & \stackrel{\phi}{\rightarrow} & {\mathcal M}
   \end{array}
\]
with $T$ an isometric $A$-$B$-bilinear embedding and $\phi$ a
bounded $A$-$B$-bilinear map. In this diagram we can replace
$\phi$ by $\phi / \| \phi \|$, a contractive map. Then there
exists a bounded $A$-linear map $\psi : \mathcal L \to \mathcal M$
such that $(\phi / \| \phi \|) = \psi \circ T$ by
\cite[Th.~2.2]{Lin1}. Since $\phi$ and $T$ are $B$-linear the map
$\psi$ turns out to be $B$-linear, too. Multiplying both sides by
the constant $\| \phi \|$ we obtain the map $\| \phi \| \psi$ that
completes the diagram above to a commutative one. So $\mathcal M$
is injective in the selected category.
\end{proof}

When
the C*-algebra of coefficients of a Hilbert C*-module $\mathcal E$
is not a unital C*-algebra and the Hilbert C*-module $\mathcal E$
is full, i.e.~its C*-algebra of coefficients $A$ is the minimal
admissible one, then we can consider $\mathcal E$ as a Hilbert
C*-module over larger C*-algebras, reasonably over C*-algebras
containing the C*-algebra of coefficients $A$ as an ideal and
belonging to the multiplier algebra $M(A)$ of $A$. However, a
construction by D.~Baki\'c and B.~Gulja{\v{s}} in \cite{BaGu1}
gives us the opportunity to establish a necessary condition on
those Hilbert $A$-modules to be injective in the category of
Hilbert $M(A)$-modules.

Let $A$ be a (non-unital) C*-algebra and $\mathcal M$ be a full
Hilbert $A$-module equipped with an $A$-valued inner product
$\langle .,. \rangle$. If $A$ is equipped with the standard
$A$-valued inner product defined by the rule $\langle a,b
\rangle_A=ab^*$, then the Hilbert $M(A)$-module ${\rm
End}^*_A(A,{\mathcal M})$ of all adjointable bounded $A$-linear
maps from $A$ to $\mathcal M$ is denoted by ${\mathcal M}_d$. The
$M(A)$-valued inner product on ${\mathcal M}_d$ is defined by
$\langle r,s \rangle = s^* \circ r$ for any $r,s \in {\mathcal
M}_d$. One of the remarkable properties of this construction is
the existence of an isometric embedding $\Gamma$ of $\mathcal M$
into ${\mathcal M}_d$. It is defined by the
formula $\Gamma(x)(a) = ax$ for any $a \in A$, each $x \in
\mathcal M$. The image $\Gamma({\mathcal M}) \subseteq {\mathcal
M}_d$ coincides with the subset $A \cdot {\mathcal M}_d$. Note,
that the construction depends on the unitary equivalence classes
of both the $A$-valued inner products on $A$ and on $\mathcal M$.
Furthermore, ${\mathcal M}_d$ can be characterized topologically
as the linear hull of the completion of the unit ball of $\mathcal M$
with respect to the strict topology, where the strict topology is
induced by the set of semi-norms $\{ \| \langle \cdot,x \rangle
\|_A : x \in {\mathcal M} \} \cup \{ \| b \cdot \|_{\mathcal M} :
b \in A \}$. So the described extension turns out to be a closure
operation, i.e., ${\mathcal M}_d \equiv {{\mathcal M}_d}_d$ for
any Hilbert C*-module $\mathcal M$. Finally, the closure operation
obeys orthogonal decompositions, i.e. $({\mathcal M} \oplus
{\mathcal N})_d = {\mathcal M}_d \oplus {\mathcal N}_d$, and the
sets of all adjointable bounded module maps on $\mathcal M$ and on
${\mathcal M}_d$ are always $*$-isomorphic, simply by restricting
operators on ${\mathcal M}_d$ to the $M(A)$-invariant subset
$\Gamma ({\mathcal M}) \subseteq {\mathcal M}_d$ that is
isometrically isomorphic to $\mathcal M$. For all these results we
refer to \cite{BaGu1}.

\begin{proposition}
   Let $A$ be a non-unital C*-algebra and $\mathcal E$ be a full
   Hilbert $A$-module. Let $B$ be another C*-algebra that admits a
   $*$-repre\-sen\-tation on $\mathcal E$. If $\mathcal E$ is
   injective in $\mathcal B(M(A),B)$, then ${\mathcal E} \equiv
   {\mathcal E}_d$.
\end{proposition}

\begin{proof}
Note, that the isomorphism of the sets of all adjointable bounded
module maps on both the Hilbert $A$-module and on its strict
closure turns the strict closure into a $M(A)$-$B$ bimodule, too.
So $\mathcal E_d$ is contained in the same category under
consideration.

Referring to the definition of injectivity, set $\mathcal M =
\mathcal E$, ${\mathcal N} = {\mathcal E}_d$, identify $\mathcal E$
with its image, $\Gamma(\mathcal E) \subseteq \mathcal E_d$ and $\phi = {\rm
id}_{\mathcal E}$. Since $\Gamma(\mathcal E)$ is injective, there is a
bounded $M(A)$-$B$-bimodule map, $\psi: \mathcal E_d \to \Gamma(\mathcal E)$
extending the identity map. Furthermore, by \cite[Th.~6.4]{Fr98}
we have the canonical isometric inclusions $\mathcal E \hookrightarrow
{\mathcal E}_d \hookrightarrow {\mathcal E}'$, and the $M(A)$-linear bounded
identity operator on $\mathcal E$ has a unique extension to the
identity operator on ${\mathcal E}'$ preserving the norm. In
particular, the identity operator on $\mathcal E$ extends uniquely to
the identity operator on ${\mathcal E}_d$. Therefore, ${\mathcal E}
\equiv {\mathcal E}_d$.
\end{proof}

We have that ${\mathcal E} \equiv {\mathcal E}_d$ for
a Hilbert $A$-module $\mathcal E$ provided that either the
C*-algebra $A$ of coefficients or the C*-algebra ${\rm
K}_A({\mathcal E})$ is unital. Whether only
direct orthogonal sums of Hilbert $A$-modules of these two types
can possess this closure property is an open problem at present,
cf.~\cite{Ba3}.

\begin{corollary}
   Let $A$ be a C*-algebra. If $A$ is injective in the category
   $\mathcal B(M(A), \mathbb C)$, then $A$ has to be unital (i.e. $A = M(A)$) and
   monotone complete. Moreover, if $A^N$ is injective for some
   $N \in \mathbb N$, then $A^N$ is injective for any $N \in
   \mathbb N$, in particular for $N=1$.
\end{corollary}

\begin{proof}
This follows from the facts that $A_d = M(A)$ and $(A^N)_d =
M(A)^N$ for any $N \in \mathbb N$ by construction. So $A = M(A)$
by the previous proposition. Furthermore, Proposition
\ref{prop-example} and Theorem \ref{thm-criterion} force $A$ to be
monotone complete.
\end{proof}

\section{Projectivity}

Let $A$ and $B$ be two fixed C*-algebras. We consider the
categories $\mathcal B(A,B)$ (respectively, $\mathcal B^*(A,B)$) consisting
of Hilbert $A$-$B$-bimodules as objects and
(adjointable) bounded $A$-$B$-bilinear maps as morphisms, where
$A$ serves as the C*-algebra of coefficients and $B$ admits a
$*$-representation in the C*-algebras of all adjointable bounded
operators on Hilbert $A$-modules that are objects. By definition a
Hilbert $A$-$B$-bimodule $\mathcal F$ is {\it projective} in
$\mathcal B(A,B)$ (respectively, $\mathcal B^*(A,B)$) if and only if the diagram
\begin{equation} \label{diag-proj}
   \begin{array}{lcc}
       {\mathcal N} & & \\
       \downarrow\vcenter{\rlap{$T$}} & & \\
       {\mathcal M} & \stackrel{\phi}{\leftarrow} & {\mathcal F}
   \end{array}
\end{equation}
can be completed to a commutative one by an $A$-$B$-(respectively,
adjointable) bimodule morphism $\psi: \mathcal F \to \mathcal N$,
whenever $T$ is a surjective $A$-$B$-bimodule (respectively, adjointable)
morphism and $\phi$ is a (respectively, adjointable) $A$-$B$-bimodule
morphism between Hilbert $A$-$B$-bimodules.

It is fairly easy to prove (and we do) that every object is
projective in $\mathcal B^*(A,B)$. We do not know if the same is true for
$\mathcal B(A,B)$, but we identify a family of C*-algebras for which every
object in $\mathcal B(A,B)$ is projective.

We begin by disposing of the $\mathcal B^*(A,B)$ case.

\begin{theorem} \label{basic-proj}
   Let $A$ be an arbitrary C*-algebra and $\{ {\mathcal F}, \langle .,.
   \rangle \}$ be a Hilbert $A$-module.  Let $B$ be another C*-algebra that
   admits a $*$-repre\-sen\-tation in ${\rm End}_A^*({\mathcal F})$.
   Then $\mathcal F$ is a projective object in the category $\mathcal B^*(A,B)$.
\end{theorem}

\begin{proof}
Consider an adjointable surjective bounded $A$-$B$-bimodule map
$T:\mathcal N \to \mathcal M$ of two Hilbert $A$-$B$-bimodules
$\mathcal M$ and $\mathcal N$. Since $T$ possesses closed range by
definition, the range of $T^*:\mathcal M \to \mathcal N$ is closed
in $\mathcal N$ and an orthogonal summand by Proposition
\ref{closed-range}. Since $T$ is surjective, $T^*$ has to be
injective, and we have the decomposition \linebreak[4] ${\mathcal N} =
T^*(\mathcal M) \oplus {\rm Ker}(T)$. Both these orthogonal
summands are $A$-$B$-invariant by construction. Every element $x
\in \mathcal M$ possesses a unique pre-image $T^{-1}(x) \in
T^*(\mathcal M)$. The operator $T^{-1}: \mathcal M \to
T^*(\mathcal M) \subseteq \mathcal N$ defined this way is
everywhere defined on $\mathcal M$ and possesses a closed range,
hence, it is bounded. Moreover, it is $A$-$B$-linear. Setting
$\psi: \mathcal F \to \mathcal N$ to be defined by the rule
$\psi(f)=T^{-1}(\phi(f)) \in T^*(\mathcal M) \subseteq \mathcal N$
for $f \in \mathcal F$ we obtain a bounded $A$-$B$-bilinear map
$\psi$ completing the diagram (\ref{diag-proj}) to a commutative
one.
\end{proof}

The following test for projectivity is often useful.

\begin{theorem}   \label{lem-topsumm2}
   Let $A$ and $B$ be arbitrary C*-algebras and $\{ \mathcal F,
   \langle .,. \rangle \}$ be a Hilbert $A$-$B$-bimodule. Then the
   following statements are equivalent:
   \begin{itemize}
   \item[(i)] $\mathcal F$ is projective in $\mathcal B(A,B)$,
   \item[(ii)] every bounded, surjective bimodule map, $T: \mathcal N \to
     \mathcal F$ has a right inverse, $S: \mathcal F \to \mathcal M$ that is a bounded
     bimodule map,
   \item[(iii)] whenever, $T: \mathcal N \to \mathcal F$ is a bounded, surjective
     bimodule map, $Ker(T)$ is a topological summand of $\mathcal N$ with a
     complementary space that is a bimodule, i.e.~is a topological
     bimodule summand.
   \end{itemize}
\end{theorem}

\begin{proof} The equivalence of (ii) and (iii) is clear.

Assume that $\mathcal F$ is projective.
By definition there exists an $A$-$B$-bimodule morphism $\psi:
\mathcal F \to \mathcal N$ such that $T \circ
\psi  = {\rm id}_{\mathcal F}$. By \cite[Lemma 3.1.8(2)]{Kasch} we
have the set identities $\mathcal N = T^{-1}(\mathcal F) = {\rm
Im}(\psi) + {\rm Ker}(T)$ and \linebreak[4] $\{ 0 \} = \psi({\rm
Ker} ({\rm id}_{\mathcal F})) = {\rm Im}(\psi) \cap {\rm Ker}(T)$.
Therefore, the Hilbert $A$-$B$-bimodule ${\rm Ker}(T) \subseteq
\mathcal N$ is a topological summand with topological complement
${\rm Im}(\psi)$ there, i.e.~$\mathcal N = \psi(\mathcal F)
\stackrel{.}{+} {\rm Ker}(T)$. The invariance of ${\rm Ker}(T)$
under the action of $B$ is caused by the $A$-$B$-bilinearity of
the operator $T$. Thus, (i) implies (iii).

Conversely, assume that (ii) holds and consider the situation of
diagram (1). Let $\mathcal L= \{ (f,n) \in \mathcal F \oplus \mathcal N: \phi(f)=
T(n) \}$, which is an $A$-$B$-submodule of $\mathcal F \oplus \mathcal N$.
The map $R: \mathcal L \to \mathcal F$, defined by $R((f,n))= f$ is a bounded
bimodule surjection and hence has a right inverse, $S: \mathcal F \to
\mathcal L$. Let $P: \mathcal L \to \mathcal N$ be defined by $P((f,n))= n$, so $P$
is also a bounded bimodule map and $\psi = P \circ S: \mathcal S \to
\mathcal N$ is the desired lifting of $\phi$.
\end{proof}

Theorem \ref{lem-topsumm2} indicates a way to find non-projective
Hilbert $A$-modules if such Hilbert C*-modules exist at all.

\begin{theorem} \label{c0-proj}
   Let $A$ be a C*-algebra of type $c_0$-$\sum_i \oplus {\rm K}(H_i)$,
   i.e.~a C*-algebra of compact operators on a certain Hilbert space.
   Let $\{ \mathcal F, \langle .,. \rangle \}$ be a Hilbert $A$-module
   and $B$ be another C*-algebra admitting a $*$-representation in
   ${\rm End}_A^*(\mathcal F)$. Then $\mathcal F$ is a projective
   object in $\mathcal B(A,B)$.
\end{theorem}

\begin{proof}
Referring to Theorem \ref{Mag-Schw} and Proposition \ref{Fr-neu}
we see that every kernel is an orthogonal summand. When the map is a
bimodule map, the kernel and its orthogonal complement are both bimodules.
\end{proof}

\begin{problem}
 Are the C*-algebras $A$ of type $c_0$-$\sum_i {\rm K}(H_i)$ the only
 C*-algebras $A$ for which all Hilbert $A$-modules are projective in
 $\mathcal B(A,\mathbb C)$, or not?
\end{problem}

In fact, we do not even know whether or not every Hilbert
$A$-$B$-bimodule is projective in $\mathcal B(A,B)$ for every pair of
C*-algebras, $A$ and $B$. The investigations of the authors
did not reveal any counterexample, so we state the question as a
problem to the readers:

\begin{problem}
 Does there exist a C*-algebra for which there is a non-projective
 Hilbert $A$-module in the category $\mathcal B(A,\mathbb C)$? Does there
 exist a pair of C*-algebras $A$, $B$ and a non-projective Hilbert
 $A$-$B$-bimodule in $\mathcal B(A,B)$?
\end{problem}

By \ref{lem-topsumm2}, the above problem is equivalent to
determining whether or not every surjective  bimodule map
between Hilbert bimodules has a right inverse that is a bimodule map.
\newline
The following general result partially links the final
solution of the projectivity problem to the solution of Problem
\ref{problemYY} above:

\begin{corollary}
   Let $A$ be a C*-algebra.
   Every Hilbert $A$-module is projective in the category
   $\mathcal B(A, \mathbb C)$ if and only if the kernel of every surjective
   bounded $A$-linear map between Hilbert $A$-modules is a
   topological summand.
\end{corollary}

\begin{proof}
Apply \ref{lem-topsumm2}(iii).
\end{proof}

We now take a closer look at projectivity in the case of
unital C*-algebras and its connection with Kasparov's stabilization
theorem.

\begin{proposition} \label{fam-proj}
   Let $A$ be a unital C*-algebra. Then for every $N \in \mathbb N$
   the Hilbert $A$-module, $A^N$ is projective in $\mathcal B(A, \mathbb C)$.
\end{proposition}

\begin{proof}
Given a Hilbert $A$-module $\mathcal N$ and a bounded surjective
module map, $T: \mathcal N \to A^N$, choose elements, $x_j \in \mathcal N$,
such that $T(x_j)=e_j$, where $e_j$ denotes the element that
is $1_A$ in the $j$-th component and 0, elsewhere. The map
$R: A^N \to \mathcal N$, defined by $R((a_1,...,a_N))= \sum_j a_jx_j$
is a right inverse for $T$.
\end{proof}

\begin{problem}
   When is a non-unital C*-algebra, $A$, a projective object in
   $\mathcal B(A, \mathbb C)$?
\end{problem}

By \ref{c0-proj}, some non-unital C*-algebras are projective in
$\mathcal B(A, \mathbb C)$. Also, by the above result it is easy to see
that any time $A$ is projective, then $A^N$ is projective.

The corresponding infinite dimensional version of $A^N$ is
$\ell^2(A)= \{ (a_1,a_2,...): \sum_{n=1}^{\infty} a_na_n^* \in A
\}$, where the convergence is in the norm sense.

\begin{proposition}
   If $\ell^2(A)$ is projective in $\mathcal B(A, \mathbb C)$, then
   every countably generated Hilbert $A$-module is projective
   in $\mathcal B(A, \mathbb C)$.
\end{proposition}

\begin{proof}
If $\mathcal M$ is countably generated then by Kasparov's stabilization
theorem \cite{Kas}, $\mathcal M \oplus \ell^2(A)$ is $A$-module
isomorphic to $\ell^2(A)$. Thus, $\mathcal M$ is isomorphic to an
orthogonally complemented submodule of $\ell^2(A)$. Now an elementary
diagram chase shows that an orthogonally complemented submodule of
a projective module is projective.
\end{proof}

\begin{problem}
   Let $A$ be a C*-algebra, when is $\ell^2(A)$ projective in
   $\mathcal B(A, \mathbb C)$?
\end{problem}

We will make some progress on this question below. For these results
we will need some concepts from operator spaces. Given any Hilbert
C*-module, $\mathcal M$ we can represent it as operators on a Hilbert
space. This allows us to make sense of the norms of matrices over
the Hilbert C*-module and these norms turn out to be canonical,
i.e.~to only depend on the inner product. For our purposes, we
will only need to refer to $M_{\infty}(A)$ which denotes the set
of $\infty \times \infty$ matrices over $A$ which are bounded,
i.e.~such that $\|(a_{i,j})\| \equiv \sup_{n} \|(a_{i,j})_{i,j=1}^n
\| < +\infty$ and to $C_{\infty}(\mathcal M)= \{ (m_1, m_2,....)^t:
( \langle m_i,m_j \rangle) \in M_{\infty}(A) \}$.

\begin{proposition}
   Let $\phi: \ell^2(A) \to \mathcal M$ be defined by $\phi((a_1,a_2,...))
   = \sum_n a_nm_n$. Then $\phi$ defines a bounded $A$-module map if
   and only if $\|( \langle m_i,m_j \rangle )\|$ is finite. Moreover,
in this case, $\|\phi\|= \|(\langle m_i,m_j \rangle )\|$.
\end{proposition}

\begin{proof}
For any finitely supported tuple, we have
$\|\phi((a_1,...,a_n,0,0...))\| = \| \sum_{i,j=1}^n a_i \langle
m_i,m_j \rangle a_j^* \|$.  But for any $(p_{i,j}) \in M_n(A)$, we
have that
\[
   \|(p_{i,j})\| = sup \{ \| \sum_{i,j=1}^n a_ip_{i,j}a_j^*
   \|: \sum_{j=1}^n a_ja_j^* \le 1_A \} \, ,
\]
from which the result follows.
\end{proof}

\begin{theorem}
   Let $A$ be a unital C*-algebra. Then $\ell^2(A)$ is projective
   in $\mathcal B(A,\mathbb C)$ if and only if for every pair of Hilbert
   $A$-modules, $\mathcal N, \mathcal M$ and every bounded, surjective module map,
   $T: \mathcal N \to \mathcal M$, the induced map $T_{\infty}: C_{\infty}(\mathcal N)
   \to C_{\infty}(\mathcal M)$, is surjective.
\end{theorem}

\begin{proof}
Assume that we are in the setting of diagram (1). Since the map
$\phi: \ell^2(A) \to \mathcal M$ is bounded, we have $(m_1,m_2,...)^t
\in C_{\infty}(\mathcal M)$, with $\phi((a_1,...)) = a_1m_1 + ...$, and
in order to lift $\phi$ to a map $\psi$ we must find $(n_1, ...)^t
\in C_{\infty}(\mathcal N)$, with $T(n_i)=m_i$, for all $i$.
\end{proof}

Note that we do not require the map $T_{\infty}$ to be
bounded in the above result, only onto.

We now take a closer look at what projectivity implies
for non-unital C*-algebras.

\begin{corollary}
  Let $A$ be a non-unital C*-algebra. If $A$ equipped with the
  canonical $A$-valued inner product is a projective Hilbert
  $A$-module in the category $\mathcal B(A,\mathbb C)$, then every element
  $t \in {\rm LM}(A)$ that induces a surjective map $T:A \to A$
  by the formula $T(a)=at^*$, admits a right inverse that is an
  element of ${\rm LM}(A)$, and the kernel of $T$ is a topological
  summand of $A$. Moreover, every surjective bounded module map
  $T:A \to A$ is realized by multiplication by a left multiplier
  in the way indicated.
  \newline
  If ${\rm M}(A)={\rm LM}(A)$ for the C*-algebra under consideration,
  then these conditions are automatically fulfilled.
\end{corollary}

\begin{proof}
Consider the diagram (\ref{diag-proj}) setting ${\mathcal
N}={\mathcal M}= A$ and $\phi=id_A$. Since $A$ is supposed to be a
projective Hilbert $A$-module there exists a map $\psi: A \to A$
which is implemented by the rule $\psi(a)=as^*$ for some $s \in
{\rm LM}(A)$ by the existing canonical identification of ${\rm
End}_A(A)$ with ${\rm LM}(A)$, cf.~\cite{Lin:90/2}. Note, that $T
\circ \psi = \phi$ by the choice of $\psi$. Consequently, $1_A =
1_{LM(A)}= s^*t^*=ts$ since $a=1_A$ is a possible choice for the
free variable. So $stst=s(ts)t=st$ and the element $p=st$ is an
idempotent element of ${\rm LM}(A)$. So $s \in {\rm LM}(A)$ is the
right inverse of $t \in {\rm LM}(A)$. Note, that the
idempotent $(1_A-p) \in {\rm LM}(A)$ maps $A$ onto the kernel of
the map $T$ which becomes a topological summand of the Hilbert
$A$-module $A$. The last two statements follow from the canonical
identification of ${\rm End}_A(A)$ with ${\rm LM}(A)$ and from
spectral decomposition in ${\rm M}(A)$, cf.~\cite{Lin:90/2} and
Proposition \ref{closed-range}.
\end{proof}



We close by looking at what projectivity means in the purely
algebraic category consisting of all $A$-modules and of all
$A$-linear maps for respective Hilbert $A$-modules. We show that
finitely generated Hilbert $A$-modules are in fact also projective
in all the categories of Hilbert $A$-$B$ bimodules under
consideration, as one might expect.

\begin{theorem}
   Let $A$ be a unital C*-algebra and $\mathcal F$ be a finitely
   generated Hilbert $A$-module (which is automatically a
   projective object in the category consisting of all $A$-modules
   over a fixed C*-algebra $A$ and of all $A$-linear maps by
   \cite[15.4.8]{NEWO}). Let $B$ be a C*-algebra represented as a
   C*-algebra of bounded adjointable operators on $\mathcal F$.
   Then $\mathcal F$ is an orthogonal summand of some Hilbert
   $A$-module $A^n$, $n < \infty$, and $\mathcal F$ is projective
   in the categories $\mathcal B(A,B)$ and in $\mathcal B^*(A,B)$.
\end{theorem}

\begin{proof}
Fix an $A$-valued inner product $\langle .,. \rangle$ on $\mathcal
F$. By \cite[Cor.~15.4.8]{NEWO} and by the definition of
projective $A$-modules in algebra $\mathcal F$ has to be finitely
generated, and every finitely generated Hilbert $A$-module is
projective in the purely algebraic sense. Consider the diagram
(\ref{diag-proj}) again. By supposition there exists an $A$-linear
map $\psi: \mathcal F \to \mathcal N$ such that $\phi = T \circ
\psi$. Let us show that $\psi$ is bounded. By \cite{LaFr98} there
exists a finite algebraic set of generators $\{ x_1, ... ,x_n \}$
of $\mathcal F$ such that the reconstruction formula $x =
\sum_{i=1}^n \langle x,x_i \rangle x_i$ is valid for every $x \in
\mathcal F$. This set $\{ x_1, ... ,x_n \}$ of generators is
called a normalized tight frame of $\mathcal F$ with respect to
the fixed $A$-valued inner product $\langle .,. \rangle$.
Therefore, $\psi(x) = \sum_{i=1}^n \langle x,x_i \rangle
\psi(x_i)$ for any $x \in \mathcal F$. Using the Cauchy-Schwarz
inequality for Hilbert C*-modules (\cite[Prop.~1.1]{Lance}) we
obtain the inequality

\begin{eqnarray*}
   \| \psi(x) \| & = & \left\| \left\langle
         \sum_{i=1}^n \langle x,x_i \rangle \psi(x_i) , \psi(x)
         \right\rangle_{\mathcal N} \right\| \\
         & = & \left\|
         \sum_{i=1}^n \langle x,x_i \rangle \langle \psi(x_i) , \psi(x)
         \rangle_{\mathcal N} \right\| \\
         & \leq &
         \sum_{i=1}^n \|x\|^{1/2} \| x_i \|^{1/2} \|
         \psi(x_i) \|_{\mathcal N}^{1/2} \| \psi(x) \|_{\mathcal N}^{1/2} \\
         & = &  \left(
         \sum_{i=1}^n  \| x_i \|^{1/2} \| \psi(x_i) \|_{\mathcal N}^{1/2}
         \right) \|x\|^{1/2} \| \psi(x) \|_{\mathcal N}^{1/2} \, .\\
\end{eqnarray*}
Cutting by $\| \psi(x) \|_{\mathcal N}^{1/2}$ the boundedness of
$\psi$ and, hence, the assertion of the theorem becomes obvious.
\end{proof}

\begin{problem}     {\rm
  Prove or disprove that selfdual Hilbert C*-modules are projective
  objects in the categories under consideration.  }
\end{problem}


\smallskip
The problem of finding non-projective Hilbert C*-modules in
the category $\mathcal B(A,\mathbb C)$ is closely related to the problem of characterizing
surjective bounded $A$-linear maps between Hilbert C*-modules that
do not admit right inverses in the set of all bounded
$A$-linear maps. Note, that in case the
domain and the range of the surjective maps are identified such maps
can be considered as special left multipliers of the C*-algebra of
all 'compact' operators on the underlying Hilbert C*-module, which
turns the problem into an open C*-algebraic problem of left
multiplier algebras of C*-algebras.

\smallskip
In fact, if a bounded module map $T: {\mathcal M} \to
{\mathcal N}$ is surjective then the right ideal $T \cdot {\rm
K}_A({\mathcal M} \oplus {\mathcal N})$ of the C*-algebra of
'compact' operators ${\rm K}_A({\mathcal M} \oplus {\mathcal N})$
is closed and hence,
\begin{eqnarray*}
   0 \to T \cdot {\rm K}_A({\mathcal M} \oplus {\mathcal N}) & \to &
   {\rm K}_A({\mathcal M} \oplus {\mathcal N}) \to \\
   & \to & {\rm K}_A({\mathcal M} \oplus {\mathcal N})/ T \cdot {\rm K}_A
   ({\mathcal M} \oplus {\mathcal N}) \to 0
\end{eqnarray*}
is a short exact sequence of Hilbert ${\rm K}_A({\mathcal M} \oplus
{\mathcal N})$-modules. (Here $T$ is identified with the operator
$T \oplus 0$ of ${\rm End}_A({\mathcal M} \oplus {\mathcal N})$,
and the multiplier C*-algebra of ${\rm K}_A({\mathcal M} \oplus
{\mathcal N})$ is identified with the set of all adjointable
bounded module maps on a copy of itself, cf.~\cite{Kas}.)
The sequence above would be split, i.e.~the Hilbert ${\rm
K}_A({\mathcal M} \oplus {\mathcal N})$-submodule $T \cdot {\rm
K}_A({\mathcal M} \oplus {\mathcal N})$ would be a topological
summand of the Hilbert ${\rm K}_A({\mathcal M} \oplus {\mathcal
N})$-module ${\rm K}_A({\mathcal M} \oplus {\mathcal N})$, if and
only if the operator $T$ would admit a right inverse in the
Banach algebra of all bounded module maps on ${\rm
K}_A({\mathcal M} \oplus {\mathcal N})$, if and only if the
norm-closed left ideal $T \cdot {\rm K}_A({\mathcal M} \oplus
{\mathcal N})$ can alternatively be characterized as an ideal of
the form $P \cdot {\rm K}_A({\mathcal M} \oplus {\mathcal N})$ for
some idempotent bounded module map $P$ on the Hilbert
C*-module ${\mathcal M} \oplus {\mathcal N}$, cf.~\cite{Lin:90/2}.
For the case of {\it adjointable} surjective bounded module
maps $T$ the situation is well-known: generally speaking,
adjointable bounded operators $S$ on Hilbert C*-modules $\mathcal
L$ have a norm-closed range if and only if they possess a
generalized inverse $S^+$ fulfilling $SS^+ S = S$, $S^+ S S^+ =
S^+$, if and only if the right ideal $S \cdot {\rm K}_A({\mathcal
L})$ is norm-closed, \cite{LZhang,LZhang2}. So a surjective
operator $T$ admits a generalized inverse in the C*-algebra of all
bounded modular operators on the Hilbert $K_A({\mathcal M} \oplus
{\mathcal N})$-module since the image of $T$, the set $\{ 0 \}
\oplus {\mathcal N}$, is obviously an orthogonal summand of the
Hilbert C*-module ${\mathcal M} \oplus {\mathcal N}$. For the
proofs of these facts and for the ideas on the one-sided
multiplier situation see Lun~Chuan Zhang's publications
\cite{LZhang,LZhang2}.

From this point of view, a better understanding of the properties
of non-adjoin\-table surjective bounded modular mappings would
give us much more information on the (non-)existence of
non-projective Hilbert C*-modu\-les and of non-split short exact
sequences of Hilbert C*-modules over certain C*-algebras of
coefficients. In this direction research is continuing.

\bigskip \noindent
{\bf Acknowledgements:} The authors thank Hanfeng Li for the
correction of the proof of Proposition \ref{prop-full} and for
several helpful comments which improved the paper after its
first circulation as a preprint in Winter 2006-2007.
The first author is indebted to V.~I.~Paulsen and D.~P.~Blecher
for their hospitality and support during his one-year stay at the
University of Houston in Houston, Texas, in 1998, during which
central ideas of the present paper have been worked out.



\end{document}